\documentclass[10pt]{article}

\usepackage{amsmath}
\usepackage{amssymb}
\usepackage{amscd}
\usepackage{latexsym}
\usepackage[all]{xy}

\usepackage{theorem}
\usepackage[T1]{fontenc}
\usepackage{xspace}
\newtheorem{Df}{Definition}[section]
\newtheorem{Te}[Df]{Theorem}
\newtheorem{Po}[Df]{Proposition}
\newtheorem{Cr}[Df]{Corollary}
\newtheorem{Lm}[Df]{Lemma}
\newtheorem{Ca}[Df]{Claim}
\newtheorem{Cn}[Df]{Conjecture}

\newcommand{\Bdf}{\begin{Df}}
\newcommand{\Edf}{\end{Df}}
\newcommand{\Bte}{\begin{Te}}
\newcommand{\Ete}{\end{Te}}
\newcommand{\Bpo}{\begin{Po}}
\newcommand{\Epo}{\end{Po}}
\newcommand{\Bcr}{\begin{Cr}}
\newcommand{\Ecr}{\end{Cr}}
\newcommand{\Blm}{\begin{Lm}}
\newcommand{\Elm}{\end{Lm}}
\newcommand{\Bca}{\begin{Ca}}
\newcommand{\Eca}{\end{Ca}}
\newcommand{\Bcn}{\begin{Cn}}
\newcommand{\Ecn}{\end{Cn}}

\newcommand{\beq}{\begin{equation}\label}
\newcommand{\eeq}{\end{equation}}

\def\dim{{\mbox{dim}}}

\begin{document}

\begin{center}
  \thispagestyle{empty}
{\large\bf COMBINATORICS AND $N$-KOSZUL ALGEBRAS}
\end{center}
\vspace{0.75cm}

\begin{center}
  Roland BERGER
\footnote{Laboratoire de Math\'ematiques de l'Universit\'e de Saint-\'Etienne (LaMUSE), Facult\'e des Sciences et Techniques,
  23 rue Docteur Paul Michelon,
  42023 Saint-\'Etienne Cedex 2, France\\
Roland.Berger@univ-st-etienne.fr\\
  }

\end{center} \vspace{1cm}
\begin{center}
  {\large\it in honour of Michel Dubois-Violette}

\end{center} \vspace{1cm}

\vspace {1cm}

\begin{abstract}
The numerical Hilbert series combinatorics for quadratic Koszul algebras was extended to $N$-Koszul algebras by Dubois-Violette and Popov~\cite{dvp:plac}. In this paper, we give a striking application of this extension when the relations of the algebra are all the antisymmetric tensors of degree $N$ over given variables. Furthermore, we present a new type of Hilbert series combinatorics, called comodule Hilbert series combinatorics, and due to Hai, Kriegk and Lorenz~\cite{hkl:NMMT}. When the relations are all the antisymmetric tensors, a natural generalization of the MacMahon Master Theorem (MMT) is obtained from the comodule level, the original MMT corresponding to $N=2$ and to polynomial algebras.  

\end{abstract}

\newpage

\section{Introduction}

At the beginning of 1950's, Koszul showed a remarkable homological property of polynomial algebras, and this property was used by Priddy in 1970 for defining Koszul algebras. Besides polynomial algebras, examples of Koszul algebras come from various areas such as topology (Steenrod algebra), algebraic geometry (Grassmannians, flag varieties), quantum groups (quantum spaces, quantum matrices), non-commutative algebraic geometry (quadratic Artin-Schelter regular algebras such as Sklyanin algebras). For Koszul algebras in algebraic geometry (resp. in quantum groups), the reader may consult~\cite{eis:comalg} Exercise 17.22 (resp.~\cite{man:kalg}).

In Priddy's definition, a Koszul algebra $A$ is graded, and more precisely \emph{quadratic} i.e., the generators (resp. the relations) of $A$ are homogeneous of degree 1 (resp. 2). Manin's monography~\cite{man:kalg} is an excellent introduction to quadratic algebras, including the fundamental observation that quadratic algebras form a good category for making non-commutative algebraic geometry. A more recent and more complete reference concerning quadratic algebras is the book by Polishchuk and Positselski~\cite{pp:quad}. If in the definition of quadratic algebras the relations are homogeneous of degree $N$ ($N\geq 2$), we get the class of $N$-homogeneous algebras~\cite{rbdvmw:homog}. So quadratic algebras coincide with 2-homogeneous algebras in the new terminology. We also say cubic algebras instead of 3-homogeneous algebras.

In~\cite{RB1}, we introduced a definition of Koszulity for $N$-homogeneous algebras (see below Definition \ref{defNK}), and the so-obtained algebras are said to be $N$-Koszul (we hope that this notation does not conflict with the notation used in~\cite{pp:quad} for quadratic algebras which are Koszul up to a given homological degree). Our motivating example was provided by the cubic Artin-Schelter regular algebras of global dimension 3. It is clear from the discussion initiating the classification work of Artin and Schelter~\cite{as:regular} that the minimal projective resolution of such a regular algebra satisfies an analog of the remarkable homological property discovered by Koszul for polynomial algebras. This analog takes into account the replacement of 2 by $N$ through certain jumps in the internal degrees. The sequence of these internal degrees is just the following : 
$$0,\, 1,\, N,\, N+1,\, 2N,\, 2N+1 \ldots$$  

In~\cite{RB1}, we gave a class of $N$-Koszul algebras for all integers $N\geq 2$, coinciding with polynomial algebras when $N=2$. This class is formed by the \emph{antisymmetrizer algebras} whose relations are all the antisymmetric tensors of degree $N$ over $n$ variables with $2 \leq N\leq n$. Antisymmetrizer algebras have probably some deep links with representation theory. Actually, symplectic reflection algebras~\cite{eg:reflec} can be extended to the non-quadratic case (i.e. to $N>2$) in such a way that a PBW theorem still holds, and the extension is precisely realized by the antisymmetrizer algebras~\cite{RBVG}. Other important examples of $N$-Koszul algebras come from theoretical physics such as the Yang-Mills algebras due to Dubois-Violette and Connes~\cite{CDV1, CDV2}, or from logic and combinatorial ring theory via Gerasimov's theorem~\cite{rb:gera}.

Koszul duality is an important feature of quadratic Koszul algebras and plays a crucial role in operad theory~\cite{gk:operad} or in representation theory~\cite{bgso:kdp}. In this paper, we emphasize Koszul duality from a combinatorial point of view. The first elementary connection between Koszul duality and combinatorics is probably the computation of the Hilbert series of a polynomial algebra as the inverse of the Hilbert series of the corresponding Grassmann algebra, and such a connection became rapidly clear for any quadratic Koszul algebra $A$ (see~\cite{pp:quad} and references therein). In some circumstances, there is a combinatorial interpretation of the Hilbert series of $A$ and the Koszul duality allows to compute immediately the combinatorial data thanks to the dual Koszul algebra $A^{!}$. This method to compute some combinatorial data by using Hilbert series of adapted graded algebras is called the \emph{numerical Hilbert series combinatorics}.

The numerical Hilbert series combinatorics holds for any $N$-Koszul algebra, as proved by Dubois-Violette and Popov~\cite{dvp:plac} (see also~\cite{kriegk:cras}). We shall explain below in Section 2 the result of Dubois-Violette and Popov. We shall illustrate their formula by the antisymmetrizer algebras, obtaining a beautiful combinatorial identity. When $N=2$, the identity is well-known, while it does not seem to be available in the literature when $N>2$. 

Surprisingly, the numerical Hilbert series combinatorics is linked to the MacMahon Master Theorem (MMT), and this link is situated at a rather sophisticated conceptual level. The MMT is a celebrated result in combinatorics, having a long history and subject to a lot of various proofs. Recently, Hai and Lorenz showed that the MMT is a consequence 
of another type of Hilbert series combinatorics~\cite{phhml:2MMT}, which is based on the Grothendieck ring of comodules over Manin's bialgebra of any 2-Koszul algebra (specialized to a polynomial algebra if we want to recover the original MMT). This new type of Hilbert series combinatorics can be called \emph{comodule Hilbert series combinatorics}. Comodule Hilbert series combinatorics reduces to numerical Hilbert series combinatorics by taking the dimension of comodules (comodules are assumed to be finite-dimensional as vector spaces). More recently, Hai, Kriegk and Lorenz showed that comodule 
Hilbert series combinatorics holds 
for any $N$-Koszul algebra (or superalgebra as well)~\cite{hkl:NMMT}, including certain explicit combinatorial applications (see also~\cite{ep:mmt}). 

We shall present in Section 3 the work of Hai, Kriegk and Lorenz, in particular their comodule identity. The comodule identity is reduced to the numerical identity of Dubois-Violette and Popov by taking the dimension. From the comodule identity, a matrix identity is derived, providing the so-called Koszul Master Theorem. These identities will be expressed explicitly in some striking situations, including the case of antisymmetrizer algebras. According to Etingof and Pak~\cite{ep:mmt}, a direct bijective proof of the so-obtained identities does not seem to be available in the literature. It is worth noticing that the conceptual approach of Hai, Kriegk and Lorenz allows to include the super case as well, by mixing even and odd variables.

\section{Numerical Hilbert series combinatorics}
The link between combinatorics and Koszul algebras is expressed by duality formulas
\begin{align*}
(a\ Hilbert\ & series\ associated\ to\ A) \times\\
& (a\ Hilbert\ series\ associated\ to\ A^!)=1
\end{align*}
where $A$ is a Koszul algebra and $A^!$ is the dual algebra.
Duality formulas are consequences of the Euler-Poincar\'e identity in a suitable category.
According to the category, this leads to the \emph{numerical} Hilbert series combinatorics or to the \emph{comodule} Hilbert series combinatorics.

Koszul algebras are based on a complex introduced by Koszul during the 1950's. The data are the following: $V$ is a $\mathbb{C}$-vector space with a basis $x_1,\ldots ,x_n$, $A=S(V)$ is the polynomial algebra in $x_1,\ldots ,x_n$, $A^!=\Lambda (V^{\ast})$ is the Grassmann algebra in $x^1,\ldots,x^n$ (the dual basis). Then the Koszul complex of $A$ is defined by 
$$\cdots \rightarrow S(V)\otimes \Lambda^{\ell} (V) \stackrel{d_{\ell}}{\longrightarrow} S(V)\otimes \Lambda^{\ell -1} (V) \rightarrow \cdots,$$
where for any $u_1, \ldots , u_k, v_1,\ldots , v_{\ell}$ in $V$,
\begin{align*}
&d_{\ell}((u_1\ldots u_k) \otimes(v_1\wedge \ldots \wedge v_{\ell}))  \\
&= \sum_{1\leq j \leq \ell} (-1)^{j+1} (u_1\ldots u_kv_j) \otimes(v_1\wedge \ldots \stackrel{\vee}{v_j} \ldots \wedge v_{\ell}).
\end{align*}

Koszul showed that the homology of his complex is 0 in any degree $\ell >0$.
So for any $m>0$ the subcomplex defined by $k+\ell =m$:
$$\cdots \rightarrow S^k(V)\otimes \Lambda^{\ell} (V) \stackrel{d_{\ell}}{\longrightarrow} S^{k+1}(V)\otimes \Lambda^{\ell -1} (V) \rightarrow \cdots$$
is exact and of finite length.
Applying the Euler-Poincar\'e identity relatively to the dimension map, we get 
$$\sum_{k+\ell =m} (-1)^k \dim S^k(V) \cdot \dim \Lambda^{\ell}(V)=0,$$
providing the well-known combinatorial identity 
\begin{equation} \label{wk}
\sum_{k+\ell =m} (-1)^k \binom{n+k-1}{k} \binom{n}{\ell}
 =0.
\end{equation}
Defining the Hilbert series 
$$H_A(t)=\sum_{k \geq 0} \dim S^k(V)\ t^k,$$
$$H_{A^!}(t)=\sum_{\ell \geq 0} \dim \Lambda^{\ell}(V^{\ast})\ t^{\ell},$$
all that is encoded in the single formula
$$H_A(t) \cdot H_{A^!}(-t) =1.$$

Notice that 
$A=S(V)=T(V)/(R)$, where $T(V)$ is the free algebra in $x_1, \ldots , x_n$ and $R$ is the subspace of $V\otimes V$ generated by $x_i\otimes x_j - x_j \otimes x_i$, $1\leq i<j \leq n$, and $A^!=\Lambda (V^{\ast})=T(V^{\ast})/(R^{\perp})$, where $R^{\perp}$ is the subspace of $V^{\ast}\otimes V^{\ast}$ orthogonal to $R$ for the natural pairing between $V$ and $V^{\ast}$. In both cases, $(R)$ and $(R^{\perp})$ denote the two-sided ideal generated by $R$ and $R^{\perp}$ respectively. 

This setting extends naturally as follows:
$V$ is any $\mathbb{C}$-vector space, $N$ is an integer $\geq 2$, $R$ is a subspace of $V^{\otimes N}$, $A=T(V)/(R)$. Then $A$ is a graded algebra, called an $N$-\emph{homogeneous algebra} over $V$. So $A^!=T(V^{\ast})/(R^{\perp})$ is an $N$-homogeneous algebra over $V^{\ast}$. Under these assumptions, the Koszul complex $K(A)$ of $A$ is defined and starts as  
$$\cdots \rightarrow A\otimes R \rightarrow A\otimes V \rightarrow A \rightarrow 0,$$
where the differential continues by increasing the internal degree of $1$ or $N-1$ alternately.
More precisely, if the jump map $\nu_N: \mathbb{N} \rightarrow \mathbb{N}$ is defined by 
$$\nu_N(2i)=Ni \ \ \hbox{and} \ \ \nu_N(2i+1)=Ni+1 \ \ \hbox{for}\ i\in \mathbb{N},$$
then $K(A)$ is the complex
$$\cdots \rightarrow A\otimes A^{!\ast}_{\nu_N(\ell)} \stackrel{d_{\ell}}{\longrightarrow} A\otimes A^{!\ast}_{\nu_N(\ell-1)} \rightarrow \cdots,$$
with natural arrows.
\Bdf \label{defNK} The $N$-homogeneous algebra $A$ is said to be Koszul if the homology of $K(A)$ is 0 in any degree $\ell >0$.
\Edf
The following theorem is due to Dubois-Violette and Popov~\cite{dvp:plac}.
\Bte \label{DVP} For any $N$-Koszul algebra $A$ such that the $\mathbb{C}$-vector space $V$ is finite-dimensional, we have
$$H_A (t) \cdot \left( \sum_{\ell \geq 0} (-1)^{\ell} \mathrm{dim} (A^!_{\nu_N(\ell)})\ t^{\nu_N(\ell)} \right) =1.$$
\Ete

It is a consequence of the Euler-Poincar\'e identity applied on each subcomplex of $K(A)$ corresponding to a fixed total degree. This formula was well-known for $N=2$ (Priddy, Backelin, Fr\"oberg, L\"ofwall).

\textbf{Warning.} For $N=2$, $A$ Koszul $\Rightarrow$ $A^!$ Koszul. 
But it is no longer true for $N>2$, e.g., for the Yang-Mills algebras (see Kriegk's thesis). Moreover we have only a part of $A^!$ in the formula.

\textbf{Explanation.} The good dual of $A$ is the Yoneda algebra $E(A)$ of $A$ and we have $E(A)_{\ell}\cong A^!_{\nu_N(\ell)}$ when $A$ is $N$-Koszul. In this case, $E(A)$ is an $A_{\infty}$-algebra having only two products, and the dual of the $A_{\infty}$-algebra $E(A)$ is $A$ (see~\cite{hl:higher}).

Let us give a combinatorial application of the duality formula 
$$H_A (t) \cdot \left( \sum_{\ell \geq 0} (-1)^{\ell} \dim (A^!_{\nu_N(\ell)})\ t^{\nu_N(\ell)} \right) =1,$$
when $A$ is the $N$-antisymmetrizer algebra in $x_1, \ldots ,x_n$:
$$A=\frac{\mathbb{C}\langle x_{1},\ldots ,x_{n}\rangle}{(\sum_{\sigma} \mathrm{sgn} (\sigma)\ x_{i_{\sigma(1)}}\ldots x_{i_{\sigma(N)}};\, 1\leq i_1<\cdots < i_N\leq n)}$$
where $2\leq N \leq n$. For $N=2$, $A$ is the polynomial algebra in $x_1, \ldots ,x_n$. The following was proved in~\cite{RB1}.

\Bpo \label{antK} 1) $A$ is Koszul,

2) the set of the \emph{admissible} monomials (i.e. monomials not containing any $N$-descent $x_{j_1}\ldots x_{j_N}$, $j_1>\cdots >j_N$) is a linear basis of $A$,

3) $\dim (A^!_m)=n^m$ if $0\leq m \leq N-1$, $\binom{n}{m}$ if $N\leq m \leq n$, $0$ if $m>n$.
\Epo

Denote by $L(n,N,k)$ the number of the admissible monomials of degree $k$. So
$$\sum_{k\geq 0} L(n,N,k)\ t^k= H_A(t),$$
and the duality formula gives the beautiful combinatorial identity
\begin{align*}
\sum_{k\geq 0} & L(n,N,k)\ t^k\\
& = \left( 1-nt + \binom{n}{N} t^N - \binom{n}{N+1} t^{N+1} + \binom{n}{2N} t^{2N} - \binom{n}{2N+1} t^{2N+1} + \cdots \right)^{-1}
\end{align*}
The RHS is the inverse of a polynomial whose degree is equal to $2q$ if $r=0$ and to $2q+1$ otherwise, where $n=qN+r$, $0\leq r<N$.
For $N=2$, we recover the well-known combinatorial identity (\ref{wk}) given at the beginning.

\textbf{Question (Pak).} Find a combinatorial proof when $N>2$. 
\section{Comodule Hilbert series combinatorics}

It is due to Hai, Kriegk and Lorenz~\cite{hkl:NMMT}. These authors provide for each $N$-Koszul superalgebra $A$ a \emph{matrix} combinatorial identity. If $A$ is the polynomial algebra in $x_1, \ldots , x_n$, then the MacMahon Master Theorem is recovered. For simplicity, we don't consider the super case. The setting is the following: $V$ is a $\mathbb{C}$-vector space with a basis $x_1,\ldots ,x_n$, $R$ is a subspace of $V^{\otimes N}$ where $N\geq 2$ is fixed, $A=A(V,R)$ is the associated $N$-homogeneous algebra, $A^!=A(V^{\ast}, R^{\perp})$ is the dual $N$-homogeneous algebra. 
The key ingredient is Manin's bialgebra $end(A)$
due to Manin if $N=2$~\cite{man:kalg}, and extended to $N>2$ in~\cite{rbdvmw:homog}.

Manin's bialgebra $end(A)$ is defined as being the $N$-homogeneous algebra
$$end(A)=A^! \bullet A = A(V^{\ast}\otimes V, R^{\perp} \otimes R)$$
endowed with the coproduct 
$\Delta :end(A)\rightarrow end(A)\otimes \, end(A)$ given by 

$$\Delta (z_i^j)= \sum_k z_i^k\otimes z_k^j\,,\  \mathrm{where}\  z_i^j:=x^j\otimes x_i.$$
Moreover $A$ is a left $end(A)$-comodule defined by $\delta : A \rightarrow end(A)\otimes A$ where  
$\delta (x_i)=\sum_j z_i^j \otimes x_j$, and $\delta$ is an algebra map.

\textbf{Warning.} If $A$ is the polynomial algebra in $x_1,\ldots ,x_n$, $end(A)$ is not commutative. In fact, we have only the following relations in each submatrix
$$\left( \begin{array}{cc}  
a & b \\
c & d 
\end{array} \right): \  ac=ca, \ bd=db, \ [a,d]=[c,b].$$
Thus some relations are \emph{missing} to get the commutative space $End(V)$.

The main idea is that Manin's bialgebra is sufficient to derive the MacMahon Master Theorem. The following fundamental observation is due to Hai, Kriegk and Lorenz~\cite{hkl:NMMT}.
\Bpo \label{fund} The Koszul complex $K(A)$ and all its homogeneous parts $K(A)_m$ (for the total degree $m$) are complexes of left $end(A)$-comodules. The complex $K(A)_m$ is as follows
$$\cdots \rightarrow A_k\otimes A^{!\ast}_{\nu_N(\ell)} \stackrel{d_{\ell}}{\longrightarrow} A_{k+\nu_N(\ell)-\nu_N(\ell-1)} \otimes A^{!\ast}_{\nu_N(\ell-1)} \rightarrow \cdots$$
where $k+\nu_N(\ell)=m$.
\Epo

Let $\mathcal{C}$ be the category of left $end(A)$-comodules which are finite-dimensional as vector spaces. Let $K_0(\mathcal{C})$ be the Grothendieck ring of $\mathcal{C}$.
Assume that $A$ is Koszul. Then the Euler-Poincar\'e identity applied to the exact complex $K(A)_m$ in $\mathcal{C}$ gives for any $m>0$
$$\sum_{k+\nu_N(\ell)=m} (-1)^{\ell}\, [A_k]\cdot [A^{!\ast}_{\nu_N(\ell)}] = 0.$$
For every object $M$ of $\mathcal{C}$, $[M]$ denotes the class of $M$ in the Grothendieck ring of $\mathcal{C}$.

Define the following Poincar\'e series in $K_0(\mathcal{C})[[t]]$ 
$$P_A (t)= \sum_{k\geq 0}\, [A_k]\ t^k,$$
$$P_{A^{!\ast},N}(-t)= \sum_{\ell \geq 0}\, (-1)^{\ell}\, [A^{!\ast}_{\nu_N(\ell)}]\ t^{\nu_N(\ell)} .$$
\Bte \label{comid} (Hai, Kriegk, Lorenz) If $A$ is $N$-Koszul, then the following comodule identity holds
$$P_A (t)\cdot P_{A^{!\ast},N}(-t) =1.$$
\Ete

In order to recover the numerical identity of Dubois-Violette and Popov, it suffices to apply $\dim : K_0(\mathcal{C})\rightarrow \mathbb{Z}$.
Now we want to go from the comodule identity to a matrix identity living in Manin's bialgebra (instead of the Grothendieck ring).
It is realized by the \emph{character} map 
$$\chi : K_0(\mathcal{C}) \rightarrow end(A), \ \ [M] \mapsto \chi_M.$$
Here $\chi_M$ is the image of $\delta_M$ under the natural composite map (where Tr is the trace)
$$Hom(M, end(A) \otimes M) \cong end(A)\otimes M \otimes M^{\ast} \stackrel{\mathrm{Id \otimes Tr}}{\longrightarrow} end(A).$$
Then $\chi$ is a ring map and we have the commutative diagram of ring maps
$$\xymatrix@=1cm{K_0(\mathcal{C})\ar[d]^{\mathrm{dim}}\ar[r]^{\chi}& end(A) \ar[d]^{\mathrm{counit}}\\
\mathbb{Z}\ar[r]^{\mathrm{can}}& \mathbb{C}}$$

Transporting the comodule identity from $K_0(\mathcal{C})[[t]]$ to $end(A)[[t]]$ by $\chi$, Hai, Kriegk and Lorenz get the following matrix identity, called the Koszul Master Theorem (KMT).

\Bte \label{KMT} If $A$ is $N$-Koszul, we have
$$\left(\sum_{k\geq 0} \chi(A_{k})\, t^k\right)\left(\sum_{\ell \geq 0} (-1)^{\ell} \chi(A^{!\ast}_{\nu_N(\ell)})\, t^{\nu_N(\ell)}\right)=1.$$
\Ete

To recover the MacMahon Master Theorem,
we use the notations introduced by Garoufalidis, L\`e, Zeilberger~\cite{GLZ}:
$x^{\mathbf{m}}=x_1^{m_1}\ldots x_n^{m_n}$ for $\mathbf{m}=(m_1, \ldots , m_n)\in \mathbb{N}^n$,
$X^{\mathbf{m}}=X_1^{m_1}\ldots X_n^{m_n}$ where $X_i=\sum_j z_i^j \otimes x_j$,
$G(\mathbf{m})$ is the coefficient of $x^{\mathbf{m}}$ in $X^{\mathbf{m}}$,
$Z$ is the $n\times n$ matrix with entries $z_i^j$. Moreover, the bosonic and fermionic sums are defined by

$$\mathrm{Bos}(Z)= \sum_{k\geq 0} \sum_{|x^{\mathbf{m}}|=k} G(\mathbf{m})\, t^k,$$

$$\mathrm{Ferm}(Z)= \sum_{0 \leq \ell\leq n} \sum_{J\subset \{1,\ldots, n\}, |J|=\ell}\, (-1)^{\ell} \det (Z_J) t^{\ell}.$$
Then the KMT for the polynomial algebra in $x_1,\ldots , x_n$ is equivalent to 
$$\mathrm{Bos}(Z)\cdot \mathrm{Ferm}(Z) =1,$$ 
from which the original MMT is deduced: for any complex $n\times n$ matrix $Z$,
$$\sum_{\mathbf{m}\in \mathbb{N}^n} G(\mathbf{m})\, t_1^{m_1}\ldots t_n^{m_n} =\left( \det (I_n-Z T) \right)^{-1}$$
where $T$ is the diagonal matrix with $t_1,\ldots , t_n$ in the diagonal.

The very general setting of the KMT allows to get effortless other examples with the same conceptual method, and certain examples are new:

1) If $A$ is the quantum space $x_jx_i=q_{ij}x_ix_j$, $1\leq i<j\leq n$, Hai, Kriegk and Lorenz recover the $q$-MT obtained by Garoufalidis, L\`e, Zeilberger~\cite{GLZ} in 2005 (for the one-parameter case) and by Konvalinka, Pak~\cite{kp:mmt} in 2006 (for the multi-parameter case). It is quite surprising that Andrews' problem (1975) asking for a $q$-analog of MMT was solved only thirty years later, while the quantum space was known since the 1980's. For a very interesting historical review on the subject of the non-commutative MacMahon master theorems, the reader may consult~\cite{kp:mmt}. 

2) When $A$ is the $N$-antisymmetrizer algebra, Hai, Kriegk and Lorenz recover the $N$-MT obtained by Etingof and Pak~\cite{ep:mmt} in 2006. The notations are the same as previously; moreover $G(i_1, \ldots , i_k)$ denotes the coefficient of $x_{i_1}\ldots x_{i_k}$ in $X_{i_1}\ldots X_{i_k}$, and
$\Lambda (n,N)_k$ denotes the set of admissible sequences $\mathbf{i}=(i_1,\ldots, i_k)$, i.e., without $N$-descents. Then the $N$-MT is the following

\begin{align*}
\sum_{\mathbf{i}\in \Lambda (n,N)_k, k\geq 0} &G(\mathbf{i})\, t_{i_1}\ldots t_{i_k}=\\
& \left( \sum_{0\leq r\leq n, r\equiv 0,1 (N)} (-1)^{r-0,1} c_r(ZT)\right)^{-1}
\end{align*}
where $\det(\lambda I_n-ZT)= \sum_{0\leq r\leq n} c_r(ZT)\lambda ^{n-r}$ is the expansion of the characteristic polynomial of the matrix $ZT$.

3) Hai, Kriegk and Lorenz obtain the superversion of 2), answering to a question by Konvalinka and Pak~\cite{kp:mmt}. Now $x_1,\ldots , x_p$ are even ($\hat{i}=0$), $x_{p+1}, \ldots , x_{p+q=n}$ are odd ($\hat{i}=1$), 
$\Lambda (p,q,N)_k$ is the set of superadmissible sequences $\mathbf{i}=(i_1,\ldots, i_k)$ and $\hat{\mathbf{i}}=\hat{i_1}+ \cdots + \hat{i_k}$. Then the super $N$-MT is the following
\begin{align*}
\sum_{\mathbf{i}\in \Lambda (p,q,N)_k, k\geq 0} &(-1)^{\hat{\mathbf{i}}} G(\mathbf{i})\, t^k=\\
& \left( \sum_{r\equiv 0,1 (N)} (-1)^{0,1} e_r(Z)\, t^r \right)^{-1}
\end{align*}
where $\mathrm{ber}(I_n+tZ)= \sum_{r\geq 0} e_r(Z) t^r$ is the expansion of the Berezinian (or superdeterminant) of the matrix $I_n+tZ$.


\begin{thebibliography}{99}

\bibitem{as:regular} M. Artin, W. F. Schelter, Graded algebras of global dimension 3, \emph{Adv.
Math.} \textbf{66} (1987), 171-216.

\bibitem{bgso:kdp} A. A. Beilinson, V. Ginzburg, W. Soergel, Koszul duality patterns in
representation theory, \emph{J. Am. Math. Soc.} \textbf{9} (1996), 473-527.

\bibitem{RB1} R. Berger, Koszulity for nonquadratic algebras.
\emph{J. Algebra} \textbf{239} (2001) 705-734.

\bibitem{rb:gera} R. Berger, Gerasimov's theorem and $N$-Koszul algebras, arXiv 0801.3383, to appear in \emph{J. London Math. Soc.}

\bibitem{rbdvmw:homog} R. Berger, M. Dubois-Violette, M. Wambst, Homogeneous algebras,  \emph{J. Algebra} \textbf{261} (2003), 172-185.

\bibitem{RBVG} R. Berger, V. Ginzburg, Higher symplectic reflection algebras and non-homogeneous $N$-Koszul property, \emph{J. Algebra} \textbf{304} (2006), 577-601.

\bibitem{CDV1} A. Connes, M. Dubois-Violette, Yang-Mills algebra, \emph{Lett. Math. Phys.} \textbf{61} (2002), 149-158.

\bibitem{CDV2} A. Connes, M. Dubois-Violette, Yang-Mills and some related algebras, math-ph/0411062.

\bibitem{dvp:plac} M. Dubois-Violette, T. Popov, Homogeneous algebras,
  statistics and combinatorics, \emph{Lett. Math. Phys.} \textbf{61} (2002), 159-170.

\bibitem{eis:comalg} D. Eisenbud, \emph{Commutative algebra with a view toward algebraic geometry}, GTM 150, Springer, 1995.

\bibitem{eg:reflec} P. Etingof, V. Ginzburg, Symplectic reflection
  algebras, Calogero-Moser space, and deformed Harish-Chandra homomorphism, \emph{Invent. math.} \textbf{147} (2002), 243-348. 

\bibitem{ep:mmt} P. Etingof, I. Pak, An algebraic extension of the MacMahon master theorem, \emph{Proc. Amer. Math. Soc.}  \textbf{136} (2008), 2279-2288.

\bibitem{GLZ} S. Garoufalidis, T. Lê, D. Zeilberger, The quantum MacMahon master theorem, \emph{Proc. Natl. Acad. Sci. USA} \textbf{103} (2006), 13928-13931.

\bibitem{gk:operad} V. Ginzburg, M. Kapranov, Koszul duality for operads, \emph{Duke Math. J.} \textbf{76} (1994), 203-272. Erratum in \emph{Duke Math. J.} \textbf{80} (1995), 293.

\bibitem{hkl:NMMT} P. H. Hai, B. Kriegk, M. Lorenz, $N$-homogeneous superalgebras, \emph{J. 
Noncommut. Geom.} \textbf{2} (2008), 1-51.

\bibitem{phhml:2MMT} P. H. Hai, M. Lorenz, Koszul algebras and the quantum MacMahon master theorem, 
\emph{Bull. London Math. Soc.} \textbf{39} (2007), 667-676.
 
\bibitem{hl:higher} J.-W. He, D.-M. Lu, Higher Koszul algebras and $A$-infinity algebras, \emph{J. Algebra} \textbf{293} (2005), 335-362. 

\bibitem{kp:mmt}M. Konvalinka, I. Pak, Non-commutative extensions of the MacMahon master theorem, \emph{Adv. Math.}  \textbf{216} (2007), 29-61.

\bibitem{kriegk:cras} B. Kriegk, Un crit\`ere num\'erique pour la propri\'et\'e de Koszul g\'en\'eralis\'ee, \emph{C. R. Acad. Sci. Paris}, Ser I \textbf{344} (2007), 545-548.

\bibitem{man:kalg} Y. I. Manin, \emph{Quantum groups and non-commutative geometry}, CRM,
Universit\'e de Montr\'eal, 1988.

\bibitem{pp:quad} A. Polishchuk, L. Positselski, \emph{Quadratic algebras}, ULS 37, American Mathematical Society, Providence, RI, 2005.  


\end{thebibliography}
\end{document}